\documentclass{amsart}
\usepackage{latexsym,amssymb,amsmath,amsthm,amscd,graphicx,esint}

\makeatletter
\@namedef{subjclassname@2010}{%
  \textup{2010} Mathematics Subject Classification}
\makeatother

\numberwithin{equation}{section}
\newtheorem{theorem}{Theorem}[section]

\newtheorem{proposition}[theorem]{Proposition}

\theoremstyle{definition}

\newtheorem{remark}[theorem]{Remark}

\theoremstyle{remark}

\newcommand{\cD}{{\mathcal D}}

\newcommand{\cF}{{\mathcal F}}

\newcommand{\cI}{{\mathcal I}}

\newcommand{\cW}{{\mathcal W}}

\newcommand{\R}{{\mathbb R}}

\newcommand{\Z}{{\mathbb Z}}

\def\sg{\sigma}

\def\0{\emptyset}
\def\6{\partial}
\def\8{\infty}

\def\l{\left}
\def\r{\right}
\def\ds{\displaystyle}
\def\ol{\overline}

\begin{document}

\title[The trilinear embedding theorem]{The trilinear embedding theorem}
\author[H.~Tanaka]{Hitoshi Tanaka}
\address{Graduate School of Mathematical Sciences, The University of Tokyo, Tokyo, 153-8914, Japan}
\email{htanaka@ms.u-tokyo.ac.jp}
\thanks{
The author is supported by 
the FMSP program at Graduate School of Mathematical Sciences, the University of Tokyo, 
and Grant-in-Aid for Scientific Research (C) (No.~23540187), 
the Japan Society for the Promotion of Science. 
}
\date{}

\begin{abstract}
Let $\sg_i$, $i=1,2,3$, denote 
positive Borel measures on $\R^n$, 
let $\cD$ denote the usual collection of dyadic cubes in $\R^n$ 
and let $K:\,\cD\to[0,\8)$ be a map. 
In this paper we give a characterization of the trilinear embedding theorem. 
That is, we give a characterization of the inequality 
$$
\sum_{Q\in\cD}
K(Q)\prod_{i=1}^3\l|\int_{Q}f_i\,d\sg_i\r|
\le C
\prod_{i=1}^3
\|f_i\|_{L^{p_i}(d\sg_i)}
$$
in terms of discrete Wolff's potential and 
Sawyer's checking condition, 
when $1<p_1,p_2,p_3<\8$ 
and $\frac1{p_1}+\frac1{p_2}+\frac1{p_3}\ge 1$.
\end{abstract}

\maketitle

\section{Introduction}\label{sec1}
The purpose of this paper is to investigate the trilinear embedding theorem. 
We first fix some notations. 
We will denote by $\cD$ the family of all dyadic cubes 
$Q=2^{-k}(m+[0,1)^n)$, 
$k\in\Z,\,m\in\Z^n$. 
Let $K:\,\cD\to[0,\8)$ be a map and 
let $\sg_i$, $i=1,2,3$, be positive Borel measures on $\R^n$. 
In this paper we give a necessary and sufficient condition 
for which the inequality 
\begin{equation}\label{1.1}
\sum_{Q\in\cD}
K(Q)\prod_{i=1}^3\l|\int_{Q}f_i\,d\sg_i\r|
\le C
\prod_{i=1}^3
\|f_i\|_{L^{p_i}(d\sg_i)},
\end{equation}
to hold when 
$1<p_1,p_2,p_3<\8$ and 
$\frac1{p_1}+\frac1{p_2}+\frac1{p_3}\ge 1$. 
By duality, 
\eqref{1.1} is equivalent to 
two-weight norm inequality for the bilinear positive operators 
$$
\|T_{K}[f_1d\sg_1,f_2d\sg_2]\|_{L^{p_3'}(d\sg_3)}
\le C
\prod_{i=1}^2
\|f_i\|_{L^{p_i}(d\sg_i)}.
$$
Here, for each $1<p<\8$, 
$p'$ denote the dual exponent of $p$, 
i.e., $p'=\frac{p}{p-1}$, and 
the bilinear positive operator 
$T_{K}[\cdot\sg_1,\cdot\sg_2]$ 
is given by 
$$
T_{K}[f_1d\sg_1,f_2d\sg_2](x)
:=
\sum_{Q\in\cD}
K(Q)\prod_{i=1}^2\l(\int_{Q}f_i\,d\sg_i\r)1_{Q}(x),
\quad x\in\R^n,
$$
where $1_{E}$ stands for the characteristic function of the set $E$. 

For the bilinear embedding theorem, 
in the case $\frac1{p_1}+\frac1{p_2}\ge 1$, 
Sergei Treil gives a simple proof of the following. 

\begin{proposition}[{\rm\cite[Theorem 2.1]{Tr}}]\label{prp1.1}
Let $K:\,\cD\to[0,\8)$ be a map and 
let $\sg_i$, $i=1,2$, be positive Borel measures on $\R^n$. 
Let $1<p_1,p_2<\8$ and 
$\frac1{p_1}+\frac1{p_2}\ge 1$.
The following statements are equivalent:

\begin{enumerate}
\item[{\rm(a)}] 
The following bilinear embedding theorem holds:
$$
\sum_{Q\in\cD}
K(Q)\prod_{i=1}^2\l|\int_{Q}f_i\,d\sg_i\r|
\le c_1
\prod_{i=1}^2
\|f_i\|_{L^{p_i}(d\sg_i)}
<\8;
$$
\item[{\rm(b)}] 
For all $Q\in\cD$, 
$$
\begin{cases}
\l(\int_{Q}\l(\sum_{Q'\subset Q}K(Q')\sg_1(Q')1_{Q'}\r)^{p_2'}\,d\sg_2\r)^{1/p_2'}
\le c_2
\sg_1(Q)^{1/p_1}
<\8,
\\
\l(\int_{Q}\l(\sum_{Q'\subset Q}K(Q')\sg_2(Q')1_{Q'}\r)^{p_1'}\,d\sg_1\r)^{1/p_1'}
\le c_2
\sg_2(Q)^{1/p_2}
<\8.
\end{cases}
$$
\end{enumerate}
\noindent
Moreover,
the least possible $c_1$ and $c_2$ are equivalent.
\end{proposition}

Proposition \ref{prp1.1} was first proved 
for $p_1=p_2=2$ in \cite{NaTrVo} 
by the Bellman function method. 
Later in \cite{LaSaUr}, 
this was proved in full generality. 
The checking condition in Proposition \ref{prp1.1} is called 
\lq\lq the Sawyer type checking condition", 
since this was first introduced by Eric T. Sawyer 
in \cite{Sa1,Sa2}. 

To describe the case 
$\frac1{p_1}+\frac1{p_2}<1$, 
we need discrete Wolff's potential. 

Let $\mu$ and $\nu$ be positive Borel measures on $\R^n$ and 
let $K:\,\cD\to[0,\8)$ be a map. 
We will denote 
by $\ol{K}_{\mu}(Q)(x)$ 
the function 
$$
\ol{K}_{\mu}(Q)(x)
:=
\frac1{\mu(Q)}\sum_{Q'\subset Q}
K(Q')\mu(Q')1_{Q'}(x),
\quad x\in Q\in\cD,
$$
and $\ol{K}_{\mu}(Q)(x)=0$ 
when $\mu(Q)=0$. 
For $p>1$, 
the discrete Wolff's potential 
$\cW^p_{K,\mu}[\nu](x)$ 
of the measure $\nu$ is defined by 
$$
\cW^p_{K,\mu}[\nu](x)
:=
\sum_{Q\in\cD}
K(Q)\mu(Q)
\l(\int_{Q}\ol{K}_{\mu}(Q)(y)\,d\nu(y)\r)^{p-1}
1_{Q}(x),
\quad x\in\R^n.
$$
The author prove the following, 
which describes the case 
$\frac1{p_1}+\frac1{p_2}<1$.

\begin{proposition}[{\rm\cite[Theorem 1.3]{Ta2}}]\label{prp1.2}
Let $K:\,\cD\to[0,\8)$ be a map and 
let $\sg_i$, $i=1,2$, be positive Borel measures on $\R^n$.
Let $1<p_1,p_2<\8$ and 
$\frac1{p_1}+\frac1{p_2}<1$.
The following statements are equivalent:

\begin{enumerate}
\item[{\rm(a)}] 
The following bilinear embedding theorem holds:
$$
\sum_{Q\in\cD}
K(Q)\prod_{i=1}^2\l|\int_{Q}f_i\,d\sg_i\r|
\le c_1
\prod_{i=1}^2
\|f_i\|_{L^{p_i}(d\sg_i)}
<\8;
$$
\item[{\rm(b)}] 
For $\frac1r+\frac1{p_1}+\frac1{p_2}=1$,
$$
\begin{cases}
\|\cW^{p_2'}_{K,\sg_2}[\sg_1]^{1/p_2'}\|_{L^r(d\sg_1)}
\le c_2<\8,
\\
\|\cW^{p_1'}_{K,\sg_1}[\sg_2]^{1/p_1'}\|_{L^r(d\sg_2)}
\le c_2<\8.
\end{cases}
$$
\end{enumerate}
\noindent
Moreover,
the least possible $c_1$ and $c_2$ are equivalent.
\end{proposition}

In his survey of the $A_2$ theorem \cite{Hy}, 
Tuomas P. Hyt\"{o}nen introduces 
another proof of Proposition \ref{prp1.1}, 
which uses the \lq\lq parallel corona" decomposition 
from the recent work of 
Lacey, Sawyer, Shen and Uriarte-Tuero 
\cite{LaSaShUr} 
on the two-weight boundedness of the Hilbert transform. 
In this paper, 
following Hyt\"{o}nen's arguments in \cite{Hy} and 
applying Propositions \ref{prp1.1} and \ref{prp1.2},
we shall establish the following theorem 
(Theorem \ref{thm1.3}). 

Let $\cI$ be the set of all permutations 
of $(1,2,3)$, i.e., 
$$
\cI
:=
\{
(1,2,3),\,
(1,3,2),\,
(2,1,3),\,
(2,3,1),\,
(3,1,2),\,
(3,2,1)
\}.
$$
Let $\mu$ be a positive borel measure on $\R^n$ 
and let $K:\,\cD\to[0,\8)$ be a map. 
For $Q\in\cD$, 
we will denote 
by $K(Q,\mu)$ 
the map 
$$
K(Q,\mu)(Q')
:=
\begin{cases}
K(Q')\mu(Q'),\quad q'\in\cD,\,Q'\subset Q,
\\
0,\quad\text{otherwise}.
\end{cases}
$$

\begin{theorem}\label{thm1.3}
Let $K:\,\cD\to[0,\8)$ be a map and 
let $\sg_i$, $i=1,2,3$, be positive Borel measures on $\R^n$. 
Let $1<p_1,p_2,p_3<\8$ and 
$\frac1{p_1}+\frac1{p_2}+\frac1{p_3}\ge 1$.
The following statements are equivalent:

\begin{enumerate}
\item[{\rm(a)}] 
The following trilinear embedding theorem holds:
$$
\sum_{Q\in\cD}
K(Q)\prod_{i=1}^3\l|\int_{Q}f_i\,d\sg_i\r|
\le c_1
\prod_{i=1}^3
\|f_i\|_{L^{p_i}(d\sg_i)}
<\8;
$$
\item[{\rm(b)}] 
For any $(a,b,c)\in\cI$, 
if $\frac1{p_a}+\frac1{p_b}\ge 1$, 
then we have, 
for all $Q\in\cD$, 
$$
\begin{cases}
\l(\int_{Q}\l(\sum_{Q'\subset Q}K(Q,\sg_c)(Q')\sg_a(Q')1_{Q'}\r)^{p_b'}\,d\sg_b\r)^{1/p_b'}
\le c_2
\sg_a(Q)^{1/p_a}\sg_c(Q)^{1/p_c}
<\8,
\\
\l(\int_{Q}\l(\sum_{Q'\subset Q}K(Q,\sg_c)(Q')\sg_b(Q')1_{Q'}\r)^{p_a'}\,d\sg_a\r)^{1/p_a'}
\le c_2
\sg_b(Q)^{1/p_b}\sg_c(Q)^{1/p_c}
<\8,
\end{cases}
$$
if $\frac1{p_a}+\frac1{p_b}<1$, 
then we have, 
for all $Q\in\cD$ and 
for $\frac1r+\frac1{p_a}+\frac1{p_b}=1$, 
$$
\begin{cases}
\|\cW^{p_b'}_{K(Q,\sg_c),\sg_b}[\sg_a]^{1/p_b'}\|_{L^r(d\sg_a)}
\le c_2
\sg_c(Q)^{1/p_c}<\8,
\\
\|\cW^{p_a'}_{K(Q,\sg_c),\sg_a}[\sg_b]^{1/p_a'}\|_{L^r(d\sg_b)}
\le c_2
\sg_c(Q)^{1/p_c}<\8.
\end{cases}
$$
\end{enumerate}
\noindent
Moreover,
the least possible $c_1$ and $c_2$ are equivalent.
\end{theorem}

In \cite{LiSu}, 
Kangwei Li and Wenchang Sun 
establish the corresponding results of Theorem \ref{thm1.3} 
for the bilinear fractional integrals 
in the case 
$$
\frac1{p_1}+\frac1{p_2},\,
\frac1{p_1}+\frac1{p_3},\,
\frac1{p_2}+\frac1{p_3}
\ge 1.
$$
They also treat the weak-type estimates. 
For the works using Wolff's potential,
we refer the readers to 
\cite{CaOr,CaOrVe1,CaOrVe2,CaOrVe3,Ta1,TaGu,TaTe}.

\begin{remark}\label{rem1.4}
To describe the case 
$1<p_1,p_2,p_3<\8$ and 
$\frac1{p_1}+\frac1{p_2}+\frac1{p_3}<1$,
probably, 
we need Wolff's potential of two-measures. 
But, we can not find it until now.
\end{remark}

The letter $C$ will be used for constants
that may change from one occurrence to another.

\section{Proof of Theorem \ref{thm1.3}}\label{sec2}
In what follows we shall prove Theorem \ref{thm1.3}. 
Let us start by proving that 
(a) implies (b). But, 
this is a direct consequence of Propositions 
\ref{prp1.1} and \ref{prp1.2}. 
So, we concentrate on proving that 
(b) implies (a). 
We follow the arguments due to T.~Hyt\"{o}nen in \cite{Hy} with some necessary modifications. 
We will use 
$\ds\fint_{Q}f\,d\mu$ 
to denote the integral average 
$\ds\frac1{\mu(Q)}\int_{Q}f\,d\mu$.

Let $Q_0\in\cD$ be taken large enough and be fixed. 
We shall estimate the quantity 
\begin{equation}\label{2.1}
\sum_{Q\subset Q_0}
K(Q)
\prod_{i=1}^3\l(\int_{Q}f_i\,d\sg_i\r),
\end{equation}
where $f_i\in L^{p_i}(d\sg_i)$ 
is nonnegative and is supported in $Q_0$. 

We define the collection of principal cubes 
$\cF_i$ for the pair $(f_i,\sg_i)$, 
$i=1,2,3$. Namely, 
$$
\cF_i:=\bigcup_{k=0}^{\8}\cF_i^k,
$$
where 
$\cF_i^0:=\{Q_0\}$,
$$
\cF_i^{k+1}
:=
\bigcup_{F\in\cF_i^k}ch_{\cF_i}(F)
$$
and $ch_{\cF_i}(F)$ is defined by 
the set of all \lq\lq maximal" dyadic cubes $Q\subset F$ such that 
$$
\fint_{Q}f_i\,d\sg_i
>
2\fint_{F}f_i\,d\sg_i.
$$

Observe that
\begin{alignat*}{2}
\lefteqn{
\sum_{F'\in ch_{\cF_i}(F)}\sg_i(F')
}\\ &\le
\l(2\fint_{F}f_i\,d\sg_i\r)^{-1}
\sum_{F'\in ch_{\cF_i}(F)}
\int_{F'}f_i\,d\sg_i
\\ &\le
\l(2\int_{F}f_i\,d\sg_i\r)^{-1}
\int_{F}f_i\,d\sg_i
=
\frac{\sg_i(F)}{2},
\end{alignat*}
which implies 
\begin{equation}\label{2.2}
\sg_i(E_{\cF_i}(F))
:=
\sg_i\l(F\setminus\bigcup_{F'\in ch_{\cF_i}(F)}F'\r)
\ge
\frac{\sg_i(F)}{2},
\end{equation}
where the sets 
$E_{\cF_i}(F)$, $F\in\cF_i$, 
are pairwise disjoint. 

We further define the stopping parents, 
for $Q\in\cD$, 
$$
\begin{cases}
\pi_{\cF_i}(Q)
:=
\min\{F\supset Q:\,F\in\cF_i\},
\\
\pi(Q)
:=
\l(\pi_{\cF_1}(Q),\pi_{\cF_2}(Q),\pi_{\cF_3}(Q)\r).
\end{cases}
$$
Then we can rewrite the series in \eqref{2.1} as follows:
\begin{equation}\label{2.3}
\sum_{Q\subset Q_0}
=
\sum_{(F_i)\in(\cF_i)}
\sum_{\substack{Q: \\ \pi(Q)=(F_i)}}.
\end{equation}
We notice the elementary fact that, 
if $P,R\in\cD$, then 
$P\cap R\in\{P,R,\0\}$. 
This fact implies, 
if $\pi(Q)=(F_i)$, then 
$$
Q\subset F_a\subset F_b\subset F_c
\text{ for some }
(a,b,c)\in\cI.
$$
Thus, by symmetry of the problem in \eqref{2.3}, 
we shall concentrate ourselves on the estimate 
$$
{\rm(i)}
:=
\sum_{H\in\cF_3}
\sum_{\substack{(F,G)\in(\cF_1,\cF_2): \\ F\subset G\subset H}}
\sum_{\substack{Q: \\ \pi(Q)=(F,G,H)}}
K(Q)\prod_{i=1}^3\l(\int_{Q}f_i\,d\sg_i\r).
$$
It follows that, 
for $H\in\cF_3$, 
\begin{alignat*}{2}
\lefteqn{
\sum_{F\subset G\subset H}
\sum_{\substack{Q: \\ \pi(Q)=(F,G,H)}}
K(Q)\prod_{i=1}^3\l(\int_{Q}f_i\,d\sg_i\r)
}\\ &\le
2\fint_{H}f_3\,d\sg_3
\sum_{F\subset G\subset H}
\sum_{\substack{Q: \\ \pi(Q)=(F,G,H)}}
K(Q)\sg_3(Q)
\prod_{i=1}^2\l(\int_{Q}f_i\,d\sg_i\r)
\\ &=
2\fint_{H}f_3\,d\sg_3
\sum_{F\subset G\subset H}
\sum_{\substack{Q: \\ \pi(Q)=(F,G,H)}}
K(H,\sg_3)(Q)
\prod_{i=1}^2\l(\int_{Q}f_i\,d\sg_i\r).
\end{alignat*}

We need two observations. 
Suppose that 
$\pi(Q)=(F,G,H)$ and 
$F\subset G\subset H$. 
Let $i=1,2$. 
If $H'\in ch_{\cF_3}(H)$ 
satisfies $H'\subset Q$, then, 
by definition of $\pi_{\cF_3}$, 
we must have 
\begin{equation}\label{2.4}
\pi_{\cF_3}\l(\pi_{\cF_i}(H')\r)=H.
\end{equation}
By this observation, we define 
$$
ch_{\cF_3}^i(H)
:=
\l\{H'\in ch_{\cF_3}(H):\,
\pi_{\cF_3}\l(\pi_{\cF_i}(H')\r)=H\r\}.
$$
We further observe that, 
when $H'\in ch_{\cF_3}^i(H)$, 
we can regard $f_i$ as a constant on $H'$ in the above integrals, 
that is, 
we can replace $f_i$ by $f_i^{H}$ in the above integrals,
where 
$$
f_i^{H}
:=
f_i1_{E_{\cF_3}(H)}
+
\sum_{H'\in ch_{\cF_3}^i(H)}
\fint_{H'}f_i\,d\sg_i1_{H'}.
$$

A little thought confirms that, 
by the assumption (b) and 
Propositions \ref{prp1.1} and \ref{prp1.2}, 
\begin{alignat*}{2}
\lefteqn{
\sum_{F\subset G\subset H}
\sum_{\substack{Q: \\ \pi(Q)=(F,G,H)}}
K(H,\sg_3)(Q)
\prod_{i=1}^2\l(\int_{Q}f_i^{H}\,d\sg_i\r)
}\\ &\le Cc_2
\sg_3(H)^{1/p_3}
\prod_{i=1}^2
\|f_i^{H}\|_{L^{p_i}(d\sg_i)}.
\end{alignat*}
Thus, we obtain 
$$
{\rm(i)}\le Cc_2
\sum_{H\in\cF_3}
\prod_{i=1}^2
\|f_i^{H}\|_{L^{p_i}(d\sg_i)}
\fint_{H}f_3\,d\sg_3
\sg_3(H)^{1/p_3}.
$$
Since 
$1<p_1,p_2,p_3<\8$ and 
$\frac1{p_1}+\frac1{p_2}+\frac1{p_3}\ge 1$, 
we can select the auxiliary parameters 
$s_i$, $i=1,2$, that satisfy 
$$
\frac1{s_1}+\frac1{s_2}+\frac1{p_3}=1
\text{ and }
1<p_i\le s_i<\8.
$$
It follows from H\"{o}lder's inequality with 
exponent $s_1$, $s_2$ and $p_3$ that 
\begin{alignat*}{2}
{\rm(i)}
&\le Cc_2
\prod_{i=1}^2
\l(
\sum_{H\in\cF_3}
\|f_i^{H}\|_{L^{p_i}(d\sg_i)}^{s_i}
\r)^{1/s_i}
\\ &\quad\times
\l(
\sum_{H\in\cF_3}
\l(\fint_{H}f_3\,d\sg_3\r)^{p_3}
\sg_3(H)
\r)^{1/p_3}
\\ &\le Cc_2
\prod_{i=1}^2
\l(
\sum_{H\in\cF_3}
\|f_i^{H}\|_{L^{p_i}(d\sg_i)}^{p_i}
\r)^{1/p_i}
\\ &\quad\times
\l(
\sum_{H\in\cF_3}
\l(\fint_{H}f_3\,d\sg_3\r)^{p_3}
\sg_3(H)
\r)^{1/p_3}
\\ &=: Cc_2
{\rm(i_1)}\times{\rm(i_2)}\times{\rm(i_3)},
\end{alignat*}
where we have used 
$\|\cdot\|_{l^{p_i}}\ge\|\cdot\|_{l^{s_i}}$.

For ${\rm(i_3)}$, 
using 
$\sg_3(H)\le 2\sg_3(E_{\cF_3}(H))$ 
(see \eqref{2.2}), 
the fact that 
$$
\fint_{H}f_3\,d\sg_3
\le
\inf_{y\in H}
M_{\sg_3}f_3(y)
$$
and the disjointness of the sets 
$E_{\cF_3}(H)$,
we have 
\begin{alignat*}{2}
{\rm(i_3)}
&\le C
\l(
\sum_{H\in\cF_3}
\int_{E_{\cF_3}(H)}(M_{\sg_3}f_3)^{p_3}\,d\sg_3
\r)^{1/p_3}
\\ &\le C
\l(\int_{Q_0}(M_{\sg_3}f_3)^{p_3}\,d\sg_3\r)^{1/p_3}
\le C
\|f_3\|_{L^{p_3}(d\sg_3)}.
\end{alignat*}
Here, 
$M_{\sg_3}$ is the dyadic Hardy-Littlewood maximal operator and 
we have used the $L^{p_3}(d\sg_3)$-boundedness of $M_{\sg_3}$. 

It remains to estimate ${\rm(i_1)}$. 
(${\rm(i_2)}$ can be estimated by the same manner.)
It follows that 
$$
\rm(i_1)^{p_1}
=
\sum_{H\in\cF_3}
\int_{E_{\cF_3}(H)}f_1^{p_1}\,d\sg_1
+
\sum_{H\in\cF_3}
\sum_{H'\in ch_{\cF_3}^1(H)}
\l(\fint_{H'}f_1\,d\sg_1\r)^{p_1}
\sg_1(H').
$$
By the pairwise disjointness of the sets $E_{\cF_3}(H)$, 
it is immediate that
$$
\sum_{H\in\cF_3}
\int_{E_{\cF_3}(H)}f_1^{p_1}\,d\sg_1
\le
\|f_1\|_{L^{p_1}(d\sg_1)}^{p_1}.
$$
For the remaining double sum, 
we use the definition of 
$ch_{\cF_3}^1(H)$ 
(see \eqref{2.4})
to reorganize:
\begin{alignat*}{2}
\lefteqn{
\sum_{H\in\cF_3}
\sum_{H'\in ch_{\cF_3}^1(H)}
\l(\fint_{H'}f_1\,d\sg_1\r)^{p_1}
\sg_1(H')
}\\ &=
\sum_{H\in\cF_3}
\sum_{\substack{F\in\cF_1: \\ \pi_{\cF_3}(F)=H}}
\sum_{\substack{H'\in ch_{\cF_3}(H): \\ \pi_{\cF_1}(H')=F}}
\l(\fint_{H'}f_1\,d\sg_1\r)^{p_1}
\sg_1(H')
\\ &\le
\sum_{H\in\cF_3}
\sum_{\substack{F\in\cF_1: \\ \pi_{\cF_3}(F)=H}}
\l(2\fint_{F}f_1\,d\sg_1\r)^{p_1}
\sg_1(F)
\\ &\le
\sum_{F\in\cF_1}
\l(2\fint_{F}f_1\,d\sg_1\r)^{p_1}
\sg_1(F)
\\ &\le C
\|M_{\sg_1}f_1\|_{L^{p_1}(d\sg_1)}^{p_1}
\le C
\|f_1\|_{L^{p_1}(d\sg_1)}^{p_1}.
\end{alignat*}
Altogether, we obtain 
$$
{\rm(i)}
\le Cc_2
\prod_{i=1}^3
\|f_i\|_{L^{p_i}(d\sg_i)}.
$$
This yields the theorem.

\end{document}